\documentclass[envcountsame]{mmnp}
\usepackage[tbtags]{amsmath}
\usepackage{amssymb}
\usepackage{graphicx}
\usepackage[latin1]{inputenc}

\def\Image{\mathop{\rm Im}}

\def\ep{\epsilon}
\def\d{\delta}

\def\R{{\mathbb R}}

\def\oH{{\overset{\circ}{H}}}
\def\oH1{{\overset{\circ}{H}\kern-.02in{}^1}}
\def\oH1{{\overset{\circ}{H}\kern-.02in{}^1}}

\def\bee{\begin{equation*}}
\def\eee{\end{equation*}}
\def\be{\begin{equation}}
\def\ee{\end{equation}}

\idline{Vol. 7, No. 2}{32}
\doi{10.1051/mmnp/20127102}

\numberwithin{equation}{section}

\begin{document}

\title{Inverse scattering problem for the Maxwell's equations}


\author{A. G. Ramm\inst{1}  \thanks{\email {ramm@math.ksu.edu}} }

\vspace{0.5cm}

\institute{\inst{1} Mathematics Department, Kansas State University,\\
Manhattan, KS 66506-2602, USA}


\abstract{Inverse scattering problem is discussed for the Maxwell's
equations. A reduction of the Maxwell's system to a new Fredholm
second-kind integral equation with a {\it scalar weakly singular
kernel} is given for electromagnetic (EM) wave scattering. This
equation allows one to derive a formula for the scattering amplitude
in which only a scalar function is present. If this function is
small (an assumption that validates a Born-type approximation), then
formulas for the solution to the inverse problem are obtained from
the scattering data: the complex permittivity $\ep'(x)$ in a bounded
region $D\subset \R^3$ is found from the scattering amplitude
$A(\beta,\alpha,k)$ known for a fixed $k=\omega\sqrt{\ep_0 \mu_0}>0$
and all $\beta,\alpha \in S^2$, where $S^2$ is the unit sphere in
$\R^3$,
 $\ep_0$ and $\mu_0$ are constant permittivity and magnetic
permeability in the exterior region $D'=\R^3 \setminus D$. The {\it
novel points} in this paper include:

i) A reduction of the inverse problem for {\it vector EM waves} to a
{\it vector integral equation with scalar kernel} without any
symmetry assumptions on the scatterer,

ii) A derivation of the {\it scalar integral equation} of the first
kind for solving the inverse scattering problem,

and

iii) Presenting formulas for solving this scalar integral equation.
The problem of solving this integral equation is an ill-posed one. A
method for a stable solution of this problem is given.}





\titlerunning{Inverse scattering problem for the Maxwell's equations}

\maketitle


\section{Introduction\label{s1}} There is a large literature on electromagnetic
wave (EM) scattering (\cite{[L]}- \cite{[Mu]}), to name a few books.
In \cite{[R476]} (see also [5]-[10] and references therein),
wave scattering theory is developed for small bodies of arbitrary
shapes. The idea of this paper is to reduce the Maxwell's system to
one Fredholm-type integral equation for the vector $E$ of eletric
field and use this equation  for solving inverse scattering problem.
Let $\ep'(x)=\ep(x)+\frac{i\sigma(x)}{\omega}$, where $\sigma (x)$
is the conductivity and $\ep(x)$ is the permittivity. It is assumed
that $\ep(x)=\ep_0$ in $D'=\R^3\setminus D$, i.e., outside of an
arbitrary large but finite domain $D$, $\sigma(x)=0$ outside $D$,
and the magnetic permeability $\mu=\mu_0$ in $\R^3$, where $\ep_0$
and $\mu_0$ are constants. We assume that $\ep'(x)$ is twice
continuously differentiable function.

The inverse scattering problem consists in finding complex
permittivity $\ep'(x)=\ep(x)+\frac{i\sigma(x)}{\omega}$ from the
scattering amplitude $A(\beta,\alpha, k)$, known at a fixed
wavenumber $k>0$ and all directions $\alpha$  of the incident plane
wave and all directions $\beta$ of the scattered waves. In Section 2
this integral equation \eqref{eq16a}  is derived. In Section 3
equation \eqref{eq16a} is used for a study of inverse scattering
problem. The original ideas in this paper include equation
\eqref{eq16a} and inversion procedure in Section 3, see equations
\eqref{eq34}-\eqref{eq39}.

The novel points in this paper include:

i) A reduction of the inverse problem for vector EM waves to a
vector integral equation with a scalar kernel without any symmetry
assumptions on the scatterer,

ii) A derivation of a scalar integral equation  for solving the
inverse scattering problem, and in presenting formulas for solving
this integral equation. The problem of solving this integral
equation is an ill-posed one. A method for a stable solution of this
problem is given.

The results are formulated in Lemmas 2.1 and 2.2 in Section 2 and
Theorems 1,2 in Section 3.

\section{EM wave scattering by a body\label{s2}}

Let $D$ be a finite body filled with a material that has complex
permittivity $\ep'(x)$ and constant magnetic permeability $\mu_0$.
The Maxwell's equations are:
\begin{equation}
\label{eq4}
\nabla\times E=i\omega\mu H, \quad \nabla\times
H=-i\omega\epsilon'(x)E \quad \text{in} \, \mathbb{R}^3,
\end{equation}
where $\omega>0$ is the frequency, $\mu=\mu_0$ is the constant in
$\R^3$ magnetic permeability, $\epsilon'(x)=\epsilon_0>0$ in $D'$,
$\ep_0$ is a constant,
$$\epsilon'(x)=\epsilon(x)+i\frac{\sigma(x)}{\omega},\quad
 \sigma(x)\geq 0,$$
$ \epsilon'(x)\not=0 \, \forall x \in \mathbb{R}^3$, $\epsilon'(x)
\in C^2(\mathbb{R}^3)$ is a twice continuously differentiable
function, $\epsilon(x)=\ep_0$ in $D'$, and the conductivity
$\sigma(x)=0$ in $D'$. From \eqref{eq4} one gets
\begin{equation}
\label{eq5}
\nabla\times\nabla\times E=K^2(x)E, \quad H=\frac{\nabla\times E}{i\omega\mu},
\end{equation}
\begin{equation}
\label{eq6}
K^2(x):=\omega^2 \epsilon'(x) \mu.
\end{equation}

We seek the solution of the equation
\begin{equation}
\label{eq7}
\nabla\times\nabla\times E=K^2(x)E
\end{equation}
satisfying the radiation condition:
\begin{equation}
\label{eq8}
E(x)=E_0(x) + v,
\end{equation}
where $E_0(x)$ is the plane wave,
\begin{equation}
\label{eq9}
E_0(x)=\mathcal{E} e^{ik\alpha\cdot x}, \quad
k=\frac{\omega}{c},
\end{equation}
 $c=\frac 1{\sqrt{\epsilon\mu}}$ is the wave velocity in the
homogeneous medium outside $D$, $\epsilon=\ep_0=const$ is the
dielectric parameter in the outside region $D'$, $\mu=\mu_0=const$
in $\R^3$, $\alpha\in S^2$ is the incident direction of the plane
wave, $\mathcal{E}\cdot\alpha=0,$ $\mathcal{E}$ is a constant
vector, and the scattered field $v$ satisfies the radiation
condition
\begin{equation}
\label{eq10}
\frac{\partial v}{\partial r}-ikv=o\big(\frac{1}{r}\big),\quad r=|x|\to\infty,
\end{equation}
uniformly in directions $\beta:=\frac{x}{r}.$

If $E$ is found, then the pair $\{E,H\}$, where
$$H=\frac{\nabla\times E}{i\omega\mu},$$ solves our scattering
problem.

Our goal is to derive a Fredholm second-kind integral equation for
$E$. We follow closely in this derivation  the presentation in
\cite{[R535]}. Let us rewrite equation \eqref{eq7} as
\begin{equation}
\label{eq11}
-\triangle E+\nabla(\nabla\cdot E)-k^2E-p(x)E=0,
\end{equation}
where
\begin{equation}
\label{eq12}
p(x):=K^2(x)-k^2, \quad p(x)=0 \, \text{in}\ D'.
\end{equation}
Note that $$\Image K^2(x)=\Image p(x)\geq 0.$$ This inequality will
be used in the proof of Claim 3 below (see formula  \eqref{eq23}).

It follows from \eqref{eq7} that
\begin{equation}
\label{eq13}
0=\nabla\cdot(K^2(x)E)=\nabla K^2(x)\cdot E+K^2(x)\nabla\cdot E.
\end{equation}
From \eqref{eq12} and \eqref{eq13} one gets
\begin{equation}
\label{eq14}
-\triangle E-k^2E-p(x)E-\nabla(q(x)\cdot E)=0,
\end{equation}
where
\begin{equation}
\label{eq15}
q(x):=\frac{\nabla K^2(x)}{K^2(x)}=\frac{\nabla \ep'(x)}{\ep'(x)}, \qquad
q(x)=0 \quad \text{in}\ D'.
\end{equation}

From \eqref{eq14} and \eqref{eq8} one gets
\begin{equation}
\label{eq16}
E=E_0+\int_Dg(x,y)\bigg{(}p(y)E(y)+\nabla_y(q(y)\cdot E(y))\bigg{)}dy,
\end{equation}
where
\begin{equation}
\label{eq17}
g(x,y):=\frac{e^{ik|x-y|}}{4\pi|x-y|}.
\end{equation}
If $q$ vanishes at the boundary of $D$, i.e., $\nabla K^2(x)=\nabla
p(x)$ vanishes at the boundary of $D$, then, after an integration by
parts in the last term of \eqref{eq16}, and taking into account that
$$\nabla_x g(x,y)=-\nabla_y g(x,y),$$ one gets an equivalent
equation:
\begin{equation}
\label{eq16a}
E=E_0+\int_Dg(x,y)p(y)E(y)dy +\nabla_x\int_Dg(x,y)q(y)\cdot
E(y)dy:=E_0+TE.
\end{equation}

Since we have assumed $K^2(x)\not=0$ and $K^2(x)\in
C^2(\mathbb{R}^3)$, it follows that equation \eqref{eq16a} is a
Fredholm equation of the second kind, because the integral operator
$T$ in \eqref{eq16a} is compact in $H^1(D)$, where $H^1(D)$ is the
usual Sobolev space.

Indeed, the operator $$Bf:=\int_Dg(x,y)f(y)dy$$ acts from $L^2(D)$
into $H^2(D)$, and the operator $$B_1E:=p(y)E+\nabla\big{(}q(y)\cdot
E\big{)}$$ acts from $H^1(D)$ into $L^2(D)$.  Thus, $T$ acts from
$H^1(D)$ into $H^2(D)$ and is, therefore, compact by the embedding
theorem.

Let us summarize in Lemma 2.1 the result we have proved:
\begin{lmm}
\label{Lemma1} If $K^2(x)\in C^2(\mathbb{R}^3),\ K^2(x)\not=0,\
K^2(x)=k^2>0$ in the exterior domain $D':=\mathbb{R}^3\setminus D$,
then the operator $T$ in \eqref{eq16a} is compact in $H^1(D)$, so
that equation \eqref{eq16} is of Fredholm type in $H^1(D).$
\end{lmm}

Let us now formulate the second auxiliary result in Lemma 2.2. After
its formulation the ideas of the proof are briefly described, The
proof consists of proving three claims, that are formulated in the
course of the proof.

\begin{lmm}
\label{Lemma2} Equation \eqref{eq16a} is uniquely solvable in
$H^1(D)$.
\end{lmm}

\begin{proof} Let us prove Lemma 2.2. It is sufficient to prove that the
homogeneous version of equation \eqref{eq16a} has only the trivial
solution. If $E$ solves the homogeneous equation \eqref{eq16a}, then
$E$ solves equation \eqref{eq7} and satisfies the radiation
condition \eqref{eq10}. The only solution to \eqref{eq7} satisfying
\eqref{eq10} is the trivial solution $E=0$.

Let us give details and prove the above claims.

{\it Claim 1}.\, {\it A solution to \eqref{eq16a} satisfies equation
\eqref{eq7}, and, consequently, \eqref{eq13}, that is,
$\nabla\cdot(K^2(x)E)=0$.}

Thus, Claim 1 states that \eqref{eq16a} is equivalent to the
original equation \eqref{eq7}, which is not at all obvious.

{\it Claim 2}.\, {\it The only solution to \eqref{eq7}, \eqref{eq10}
is $E=0$.}

{\it Proof of Claim 1}.\, If $E$ solves \eqref{eq16a} and is of the
form \eqref{eq8}, then
$$
(-\triangle-k^2)E=p(x)E+\nabla(q(x)\cdot E).
$$

Rewrite this equation as
$$
\nabla\times\nabla\times E-\nabla(\nabla\cdot E)-\nabla(q(x)\cdot
E)=K^2(x)E,
$$
or
\begin{equation}
\label{eq18}
\nabla\times\nabla\times E-\nabla\bigg{(}\frac{K^2(x)\nabla\cdot E+\nabla K^2(x)\cdot E}{K^2(x)}\bigg{)}
=K^2(x)E
\end{equation}
Denote
$$\frac{1}{K^2(x)}\nabla\cdot(K^2(x)E):=\psi(x).$$
Taking divergence of equation \eqref{eq18} one gets
\begin{equation}
\label{eq19}
-\triangle\psi-K^2(x)\psi=0 \quad \text{in}\,\quad \mathbb{R}^3.
\end{equation}
The function $\psi$ satisfies the radiation condition \eqref{eq10}.

{\it Claim 3}.\, {\it The only solution to \eqref{eq19}, which
satisfies condition \eqref{eq10}, is $\psi=0$.}

We prove Claim 3 below. Assuming that this claim is proved, we infer
that $\psi=0$, so $$\nabla\cdot(K^2(x)E)=0,$$ and equation
\eqref{eq7} holds.

Claim 1 is proved.      $\Box$

{\it Let us prove Claim 3}.\, Equation \eqref{eq19} can be written
as:
\begin{equation}
\label{eq20}
-\triangle\psi-k^2\psi-p(x)\psi=0 \quad \text{in} \quad \mathbb{R}^3,
\end{equation}
where $p(x)$ is defined in equation \eqref{eq12}, and $p(x)=0$ in
$D'$. It is known that if $k^2>0,\, p(x)\in L^2(D)$, $\Image
p(x)\geq 0$, $D$ is a bounded domain, and $\psi$ satisfies the
radiation condition \eqref{eq10} and equation \eqref{eq20}, then
$\psi=0$ (see, e.g., \cite{[R509]}).

For convenience of the reader, we sketch the proof. From
\eqref{eq20} and its complex conjugate one derives the relation
\begin{equation} \label{eq21}
\psi\triangle\bar{\psi}-\bar{\psi}\triangle\psi-2i\Image p(x)|\psi|^2=0.
\end{equation}
Integrate \eqref{eq21} over a ball $B_R$ centered at the origin of
radius $R$, and use the Green's formula to get
\begin{equation}
\label{eq22} \int_{|x|=R}\big{(}\psi\frac{\partial\bar{\psi}}{\partial
r}-\bar{\psi}\frac{\partial\psi}{\partial r}\big{)}ds -2i\int_D \Image
p(x)|\psi|^2dx=0.
\end{equation}
By the radiation condition \eqref{eq10} for $\psi$ one can rewrite
\eqref{eq22} as
\begin{equation} \label{eq23}
-2ik\lim_{R\to\infty}\int_{|x|=R}|\psi|^2ds-2i\int_D\Image p(x)|\psi|^2dx=0.
\end{equation}
Since $\Image p(x)\geq 0$ by our assumption, and $k>0$, it follows
that
\begin{equation} \label{eq24}
\lim_{R\to\infty}\int_{|x|=R}|\psi|^2dx=0.
\end{equation}
Relation \eqref{eq24} and the equation
\begin{equation}
\label{eq25} (\triangle+k^2)\psi=0 \quad |x|>
R_0,
\end{equation}
where $B_{R_0}\supset D$, imply $\psi=0$ for $|x|>R_0.$ (See
\cite[p.25]{[R190]}).

By the unique continuation principle for the solutions of the
homogeneous Schr\"odinger equation \eqref{eq20}, it folows that
$\psi=0$. Claim 3 is proved. $\Box$

{\it Let us prove Claim 2}.\, {\it This will complete the proof of
Lemma \ref{Lemma2}.}

If $E$ solves \eqref{eq7} and satisfies \eqref{eq10}, the the pair
$\{E,H\}$, where $H$ is defined by the second formula in
\eqref{eq5}, solves the homogeneous Maxwell's system \eqref{eq4},
satisfies the radiation condition, and $\{E,H\}$ is in
$H_{loc}^2(\mathbb{R}^3)$. It is known (see, e.g., \cite{[Mu]}) that
this implies $E=H=0.$ Lemma \ref{Lemma2} is proved.  \end{proof}

From formula \eqref{eq16a} assuming that the origin is inside $D$,
one gets:
\begin{equation} \label{eq26}
E=E_0+\frac{e^{ik|x|}}{|x|}\frac{1}{4\pi}\int_De^{-ik\beta\cdot y}
\big{(}p(y)E(y)+\nabla(q(y)\cdot E)\big{)}dy[1+O(\frac 1{|x|})]
\end{equation}
as $|x|\to\infty,\, \frac{x}{|x|}:=\beta$.

An alternative representation of the scattered field is given in
Section 3 in formula \eqref{eq28}.

\section{Inverse scattering for the Maxwell's equations\label{s3}}
We assume that the inhomogeneity is described by the permittivity
$$\ep'(x)=\ep(x)+i\frac {\sigma(x)}{\omega},$$ $\ep(x)=\ep_0$ in
$D':=\R^3\setminus D$, $\sigma(x)=0$ in $D'$, and $D$ is an
arbitrary large finite domain, $D\subset B_R$, where $B_R$ is the
ball of radius $R$ centered at the origin. The origin is inside $D$.
The inverse scattering problem consists of finding $\ep'(x)$ from
the knowledge of the scattered field measured  at large distances
from $D$. This scattered field is defined by the scattering
amplitude. Let us use equation \eqref{eq16a} for finding an unknown
$\ep'(x)$ given the scattering amplitude for all $\alpha, \beta$ and
a fixed $k>0$. Let us define the scattering amplitude
$A(\beta,\alpha,k)$ by the formula
\begin{equation} \label{eq27}
E-E_0= \frac {e^{ikr}}{r}A(\beta,\alpha,k)+o(\frac 1 r),\quad r:=|x|\to
\infty, \quad \beta:=\frac x r.
\end{equation}
From equations \eqref{eq16a} and \eqref{eq27} one gets:
\begin{equation} \label{eq28}
A(\beta,\alpha,k)= \frac 1 {4\pi}\int_De^{-ik\beta\cdot y}
p(y)E(y)dy +\frac {ik\beta}{4\pi}\int_De^{-ik\beta\cdot y}q(y)\cdot E(y)
dy.
\end{equation}
Denote by $[E,H]=E\times H$ the cross product of two vectors $E$ and
$H$, and by $(E,H)=E\cdot H$ their dot product. Then equation
\eqref{eq28} implies
\begin{equation} \label{eq29}
[\beta, A(\beta,\alpha,k)]=\frac 1 {4\pi}\int_De^{-ik\beta\cdot y}p(y)
[\beta,E(y)]dy.
\end{equation}
If one assumes that $|p(y)|\ll k^2$, and $E_0=\mathcal {E}
e^{ik\alpha\cdot x}$, then one may replace $E(y)$ in equation
\eqref{eq29} by $\mathcal {E}e^{ik\alpha\cdot x}$, and obtain in
this approximation, which is similar to the Born approximation in
quantum mechanics, the following relation:
\begin{equation} \label{eq30}
4\pi[\beta, A(\beta,\alpha,k)]=[\beta,\mathcal {E}]
\int_De^{ik(\alpha-\beta)\cdot y}p(y)dy.
\end{equation}
A novel and practically attractive feature of equation \eqref{eq30}
consists of a possibility to reduce the inverse problem of finding
an unknown permittivity in a Maxwell's system of vector equations to
solving a scalar equation  \eqref{eq33} below. Let us assume without
loss of generality that $|\mathcal {E}|=1$. Then $|[\beta,\mathcal
{E}]|=\sin \theta$, where $\theta\in [0,\pi]$ is the angle between
vectors $\beta$ and $\mathcal {E}$, $$([\beta,\mathcal {E}],
[\beta,\mathcal {E}])=\sin^2\theta,$$ and \eqref{eq30} implies
\begin{equation} \label{eq31}
4\pi([\beta, A(\beta,\alpha,k)], [\beta,\mathcal {E}])\sin^{-2}\theta:=f
=\int_De^{ik(\alpha-\beta)\cdot y}p(y)dy.
\end{equation}
The right-hand side of \eqref{eq31} does not depend on $\theta$.
Thus, its left-hand side, $f$, remains bounded for $\theta\in
[0,\pi]$. The function $$f=f(\beta, \alpha, k)=f(k(\alpha-\beta))$$
is known since $A(\beta,\alpha,k)$, $\beta$, $\theta$, and $\mathcal
{E}$ are known.

{\it Therefore, the inverse scattering problem of finding $\ep'(x)$
from $A(\beta,\alpha,k)$ is reduced to finding $p=p(y)$ from the
data $f(\beta, \alpha, k)$.}

If $p$ is found, then $\ep'$ can be found, because
\begin{equation} \label{eq32}
\ep'(x)=\frac {k^2+p(x)}{\omega^2 \mu},
\end{equation}
where $\mu=\mu_0$ in $\R^3$,  $k^2=\omega^2\ep_0 \mu_0$, and $k^2$,
$\omega$ and $\mu_0$  are assumed known.

Thus, the inverse problem is reduced  to solving the scalar equation
\begin{equation} \label{eq33}
\int_De^{ik(\alpha-\beta)\cdot y}p(y)dy=f(\beta, \alpha,k),
\end{equation}
where $k=\omega \sqrt{\ep_0\mu_0}>0$ is fixed, $\alpha,\beta\in S^2$
run independently of each other through the unit sphere $S^2$, and
the function $f=f(\beta, \alpha,k)$ is known on $S^2\times S^2$.
Let us summarize the result we have proved in Theorem 1.

{\bf Theorem 1.} {\it If the scattering amplitude $A(\beta,
\alpha,k)$ is known for a fixed $k>0$ and all $\alpha, \beta \in
S^2$, then one calculates the left-hand side of equation
\eqref{eq31} thus finding the scalar function $f(\beta, \alpha,k)$,
and then calculates the function $p(y)$ by solving equation
\eqref{eq33}.}

Let us now give a method for solving equation \eqref{eq33}. We
assume that the data $f$ is given with some error: the exact $f$ is
unknown and the "noisy" data $f_\d$ are given such that
$\|f-f_\d\|<\d$, where $\d>0$ is a small number. The problem is to
find $p=p(y)$ from equation \eqref{eq33}, given $f_\d$.

Such a problem has been studied in \cite{[R190]}. We formulate the
results from \cite{[R190]} without proof. The reader is referred to
the proofs to \cite{[R190]}, pp.259-274. The inversion formula from
\cite{[R190]}, p. 268, is:
\begin{equation} \label{eq34}
p_N(x)= \int_{S^2}\int_{S^2}f(k(s'-s))h_N(k(s'-s))e^{ik(s'-s)\cdot
x}ds ds',
\end{equation}
where $N$ is a  positive integer, $s'=\alpha$, $s=\beta$, $ds$ and
$ds'$ are surface elements of $S^2$,
\begin{equation} \label{eq35}
h_N(z):=|z|a_N(z)k^2(32\pi^4)^{-1},
\end{equation}
where
\begin{equation} \label{eq36}
a_N(z):=\int_{\R^3} \d_N(x)e^{-iz\cdot x}dx,
\end{equation}
and
\begin{equation} \label{eq37}
\d_N(x):=\big{(}1-\frac {|x|^2}{4R^2}\big{)}^N \big{(}\frac{N}{4\pi
R^2}\big{)}^{3/2}\big{(}\frac{\sin b- b\cos b}{b^3/3}\big{)}^{2N+3},\quad
b=\frac {2k|x|}{2N+3},
\end{equation}
where $R$ is the radius of the ball $B_R$ containing $D$. It is
proved in \cite{[R190]}, p.268, that
\begin{equation} \label{eq38}
\|p_N(x)-p(x)\|_{L^2(B_R)}\to 0 \quad as \quad N\to \infty.
\end{equation}
Suppose now that instead of the exact data $f$ the noisy data $f_\d$
are given, $\|f-f_d\|_{L^2(S^2\times S^2)}\le \d$. Then we compute
the approximate solution $p_{N\d}(x)$ by formula \eqref{eq30} with
$f_\d(k(s'-s))$ in place of $f(k(s'-s))$, and choose $N=N(\d)$
according to the recipe on p.269 in \cite{[R190]}. Then $N(\d)\to
\infty$ as $\d \to 0$. With this choice of $N(\d)$ one proves, as in
\cite{[R190]}, p. 269, that
\begin{equation} \label{eq39}
\|p_{N(\d)}(x)-p(x)\|_{L^2(B_R)}\to 0 \quad as \quad \d\to 0.
\end{equation}
These results were applied  by the author also to inversion of
incomplete tomographic data (see [11] and [12], pp. 259-264).

Let us summarize the result in Theorem 2.

{\bf Theorem 2.}  {\it A stable solution to equation \eqref{eq33} is
given by formulas \eqref{eq34}-\eqref{eq39}.}






\end{document}